\documentclass[a4paper, draftcls, onecolumn]{IEEEtran}

\usepackage[noadjust]{cite}
\usepackage{graphicx}
\usepackage[usenames,table]{xcolor}
\usepackage{tikzexternal}

\usepackage{paralist}
\usepackage{amsmath}
\usepackage{amssymb}
\usepackage{amsfonts}
\usepackage{algorithm}
\usepackage{algpseudocode}

\makeatletter
\let\MYcaption\@makecaption
\makeatother

\usepackage{caption}
\usepackage[font=footnotesize]{subcaption}

\makeatletter
\let\@makecaption\MYcaption
\makeatother

\newcommand{\ug}{$\mu$G}
\newcommand{\ugs}{$\mu$Gs}
\newcommand{\Es}{\ensuremath\varepsilon^{(s)}}
\newcommand{\MA}{\ensuremath\mathbf{A}}
\newcommand{\vE}{\ensuremath\mathbf{E}}
\newcommand{\ve}{\ensuremath\mathbf{e}}
\newcommand{\vlambda}{\ensuremath\boldsymbol{\lambda}}
\DeclareMathOperator{\diag}{diag}

\newtheorem{case}{Case}

\title{Distributed Energy Trading: The~Multiple-Microgrid~Case}

\author{David~Gregoratti,~\IEEEmembership{Member,~IEEE,} and Javier Matamoros%
\thanks{Authors are with the Centre Tecnol\`ogic de Telecomunicacions de Catalunya (CTTC),
Parc Mediterrani de la Tecnologia, 08860~Castelldefels, Barcelona~(Spain),
emails: first.last@cttc.es}%
\thanks{This work was partially supported by the Catalan Government under grant
2014 SGR 1567 and by the Spanish Government under grants TEC2011-29006-C03-01
and TEC2013-44591-P.}}

\IEEEspecialpapernotice{(This is an extended version of a paper with the same
title that appeared in the IEEE~Transactions~on~Industrial
Electronics, vol.~62, no.~4, pp.~2551--2559, Apr.~2015.)}

\begin{document}
\maketitle

\begin{abstract}
In this paper, a distributed convex optimization framework is
developed for energy trading between islanded microgrids. More
specifically, the problem consists of several islanded microgrids
that exchange energy flows by means of an arbitrary topology. Due to
scalability issues and in order to safeguard local information on
cost functions, a subgradient-based cost minimization algorithm is
proposed that converges to the optimal solution in a practical
number of iterations and with a limited communication overhead.
Furthermore, this approach allows for a very intuitive economics
interpretation that explains the algorithm iterations in terms of
``supply--demand model'' and ``market clearing.'' Numerical
results are given in terms of convergence rate of the algorithm and
attained costs for different network topologies.
\end{abstract}


\section{Introduction}
Worldwide energy demand is expected to increase steadily over the
incoming years, driven by energy demands from humans, industries and
electrical vehicles: more precisely, it is expected that the growth
will be in the order of 40\% by year 2030. This demand is fueled by
an increasingly energy-dependent lifestyle of humans, the emergence
of electrical vehicles as the major source of transportation, and
further automation of processes that will be facilitated by
machines.

In today's power grid, energy is produced in centralized and large
energy plants (macrogrid energy generation); then, the energy is
transported to the end client, often over very large distances and
through complex energy transportation meshes. Such a complex
structure has a reduced flexibility and will hardly adapt to the
demand growth, thus increasing the probability of grid instabilities
and outages. The implications are enormous as demonstrated by recent
outages in Europe and North America that have caused losses of
millions of Euros~\cite{Andersson2005}.

Given these problems at macro generation, it is of no surprise that
a lot of efforts have been put into replacing or at least
complementing macrogrid energy by means of local renewable energy
sources. In this context, microgrids are emerging as a promising
energy solution in which distributed (renewable) sources are serving
local demand \cite{Hatziargyriou2007}. When local production cannot
satisfy microgrid requests, energy is bought from the main utility.
Microgrids are envisaged to provide a number of benefits:
reliability in power delivery (e.g., by islanding), efficiency and
sustainability by increasing the penetration of renewable sources,
scalability and investment deferral, and the provision of ancillary
services. From this list, the capability of islanding
\cite{Balaguer2011,Sood2009,Katiraei2008} deserves special
attention. Islanding is one of the highlighted features of
microgrids and refers to the ability to disconnect
the microgrid loads from
the main grid and energize them exclusively via local energy
resources. Intended islanding will be executed in those situations
where the main grid cannot support the aggregated demand and/or
operators detect some major grid problem that may potentially
degenerate into an outage. In these cases, the microgrid can provide
enough energy to guarantee, at least, a basic electrical service.
The connection to the main grid will be restored as soon as the
entire system stabilizes again. Clearly, these are nontrivial
functionalities that may cause instability. In this regard, the
sequel of papers \cite{Guerrero2013,Guerrero2013_2} provide a
recent survey on decentralized control techniques for microgrids.

In order to improve the capabilities of the Smart
Grid, a typical approach is to
consider the case where several microgrids exchange energy to one
another 
even when the microgrids are islanded, that is disconnected from the main grid
\cite{Nunna13,Fathi13}. In other words, there exist energy flows within a
group of contiguous microgrids but not between the microgrids and the main grid. In
this context, the optimal power flow problem has recently attracted
considerable attention. For instance, in \cite{Ochoa2011} the authors
consider the power flow problem jointly with coordinated voltage
control. Alternatively, the work in \cite{Bruno2011} focuses on
unbalanced distribution networks and proposes a methodology to solve
three-phase  power flow problems based on a Newton-like descend
method. Due to the fact that these centralized solutions may suffer
from scalability \cite{Paudyal2011, Bruno2011} and privacy issues,
distributed approaches based on optimization tools have been
proposed in \cite{Kim1997, Baldick1999} and more recently in
\cite{Kraning2013,Dallanese2013, Erseghe2014}. In general, the
optimal power flow problem is nonconvex and, thus, an exact solution
may be too complex to compute. For this reason, suboptimal
approaches are often adopted. As an example, references
\cite{Lavaei2012,Dallanese2013} show that semidefinite relaxation
(see \cite{Luo2013} for further details on this technique) can, in
some cases, help to approximate the global optimal solution with
high precision. Alternatively, \cite{Erseghe2014}
resorts to the so-called alternating direction method of
multipliers (see \cite{Boyd_ADMM2010} for further information)
to solve the power flow problem in a distributed manner.

In this paper, conversely to the aforementioned works, we consider
an abstract model that allows us to focus on the trading process
rather than on the electrical operations of the grid. In terms of
trading, the massive spread of distributed
energy resources is expected to drive the transition from 
today's oligopolistic market to a more open and flexible one
\cite{Pagani2011}. This new picture of the market has triggered the
interest on new energy trading mechanisms \cite{YunpengWang14,
Nguyen12, Rahimiyan14, Kahrobaee13, Wang13,
Dimeas05,Nunna12,Nunna13,Ramachandran11}. For instance, the
authors in \cite{YunpengWang14} consider a scenario where a set
of geographically distributed energy storage units trade their
stored energy with other elements of the grid. The authors formulate
the problem as a noncooperative game that is shown to have at least
one Nash equilibrium point. In the context of demand response, the
work in \cite{Rahimiyan14} proposes an efficient energy
management policy to control a cluster of demands. Interestingly,
\cite{Nguyen12}
considers demand-response capabilities as an asset to  be offered within the market
and, thus, they can be traded between the different agents
(retailers, distribution system operators, aggregators, etc.). From
a more general perspective, the problem of energy trading between
microgrids (or, market agents) has been considered in
\cite{Kahrobaee13, Wang13, Dimeas05,
Nunna12,Nunna13,Ramachandran11}. Whereas these works mainly focus on
simulation studies and architectural issues, this paper attempts to
provide a comprehensive analytical solution for the energy trading
problem between microgrids that, besides its theoretical appeal, can
be distributedly implemented without the need of a central
coordinator. More specifically, our setting consists of $M$
microgrids in which:
\begin{inparaenum}[(i)] \item each microgrid has an associated
energy generation cost; \item there exists a cost imposed by the
distribution network operator for transferring energy between
adjacent microgrids; and \item each microgrid has an associated
power demand that must be satisfied.
\end{inparaenum} Under these considerations, we aim to find the
optimal amounts of energy to be exchanged by the microgrids in order
to minimize the total operational cost of the system (energy
production and transportation costs). Of course, a possible approach
would be to solve the optimization problem by means of a central
controller with \emph{global} information of the system. However,
such a centralized solution presents a number of drawbacks since
microgrids might be operated by different utilities and information
on production costs cannot be disclosed. Therefore, in order to
safeguard critical information on local cost functions and make the
system more scalable, we propose an algorithm based on dual
decomposition that iteratively solves the problem in a distributed
manner. Interestingly, each iteration of the resulting algorithm has
a straightforward interpretation in economical terms, once the new
operational variables we introduce are given the meaning of energy
prices. First, each microgrid locally computes the amounts of energy
it must produce, buy and sell to minimize its local cost according
to the current energy prices. Then, after exchanging the microgrid
bids, a regulation phase follows in which, in a distributed way, the
energy prices are adjusted according to the law of demand. This
two-step process iterates until a global agreement is reached about
prices and transferred energies.

The remainder of this paper is organized as follows. In Section
\ref{sec:system_model} we present the system model. Next, Section
\ref{sec:minimization} shows how the distributed optimization
framework provides a solution for the local subproblems and
gives an interpretation from an economical point of view.
Finally, in Sections \ref{sec:num_results} and \ref{sec:conclusions},
we present the numerical results and draw some conclusions.

\section{System Model}\label{sec:system_model}
Consider a system composed of $M$ interconnected microgrids (\ugs)
operating in islanded mode. During each scheduling interval, each
microgrid \ug-$i$ generates $E_i^{(g)}$ units of energy\footnote{For
simplicity, we assume that all energies correspond to constant power
generation/absorption/transfer over the scheduling interval.} and
consumes $E_i^{(c)}$ units of energy. Moreover, \ug-$i$ may be
allowed to sell energy $E_{i,j}$ to \ug-$j$, $j\ne i$, and to buy
energy $E_{k,i}$ from \ug-$k$, $k\ne i$. Then, energy equilibrium
within the \ug\ requires
\begin{equation}\label{eq:equilibrium}
E_i^{(g)}+\ve_i^T\MA^T\vE_i^{(b)} = E_i^{(c)}+ \ve_i^T\MA\vE_i^{(s)},
\end{equation}
where the two $M$-dimensional column vectors
\begin{align}\label{eq:en_vectors}
\vE_i^{(b)}&=\begin{bmatrix}E_{1,i}\\\vdots\\ E_{M,i}\end{bmatrix} &&\text{and}
& \vE_i^{(s)}&=\begin{bmatrix}E_{i,1}\\\vdots\\ E_{i,M}\end{bmatrix}
\end{align}
gather all the energies bought and sold, respectively, by \ug-$i$.
Also, we have introduced the adjacency matrix $\MA =
{[a_{i,j}]}_{M\times M}$: element $a_{i,j}$ is equal to one if there
exists a connection from \ug-$i$ to \ug-$j$ and zero otherwise. Note
that, generally, $\MA$ may be nonsymmetric, meaning that at least
two \ugs\ are allowed to exchange energy in one direction only. By
convention, we fix $a_{i,i}=0$. Also, $a_{i,j}=0\Rightarrow
E_{i,j}=0$~MWh for all $i,j=1,\dots,M$.

Next, let $C_i(E_i^{(g)})$ and $\gamma(E_{i,j})$ be the costs of producing
$E_i^{(g)}$ units of energy at \ug-$i$ and transferring $E_{i,j}$ units of
energy on the \ug-$i$--\ug-$j$ link, respectively.
This transferring cost function may model several
factors. For example, the distribution network operator, as the
enabler of the energy transfer between \ugs, may charge a tax
for energy transactions. In addition, this transfer cost function
may also account for line congestions by introducing soft constrains
on the maximum capacity of the line (see Section
\ref{ssec:cost_fncts} for further information).

Even though the solution below may be extended to the case where
different links have different transfer cost functions, we assume
that $\gamma(\cdot)$ is common among all the links in order to avoid
further complexity. As it will be clearer later, our approach is
quite general and only requires that all cost functions (both
production and transport) satisfy some mild convexity constraints,
summarized in Section~\ref{ssec:cost_fncts}.

Finally, we further assume that all \ugs\ agree to cooperate with one
another in order to minimize the total cost of the system. In other
words, the energy quantities exchanged by interconnected \ugs\ form
the equilibrium point of the following minimization problem:
\begin{equation}\label{eq:start}
\begin{aligned}
\mathcal{C}^*=\min_{\{E_{i,j}\}} & \sum_{i=1}^M
C_i\bigl(E_i^{(c)}+\ve_i^T\bigl(\MA\vE_i^{(s)}-\MA^T\vE_i^{(b)}\bigr)\bigr)
+ \sum_{i=1}^M\ve_i^T\MA^T\boldsymbol{\gamma}(\vE_i^{(b)})\\
\text{s.\ to}\;\; & E_{i,j}\ge 0, \forall i,j,\\
              & E_i^{(c)}+\ve_i^T\bigl(\MA\vE_i^{(s)}-\MA^T\vE_i^{(b)}\bigr)\ge
0, \forall i,
\end{aligned}
\end{equation}
where $\ve_i$ is the $i$-th column of the $M\times M$ identity
matrix and, with some abuse of notation, we wrote
$\boldsymbol{\gamma}(\vE_i^{(b)})=
{\begin{bmatrix}\gamma(E_{1,i})&\cdots&\gamma(E_{M,i})\end{bmatrix}}^T$.
Also, we used (\ref{eq:equilibrium}) to get rid of the variables
$\{E_i^{(g)}\}$. Note that problem (\ref{eq:start}) considers the
\ugs\ as parts of a common system (for instance, they are controlled
by the same operator) and the aim is at minimizing the global cost,
without focusing on the benefits/losses of each individual \ug. We
will see later on, however, that the proposed distributed and
iterative minimization algorithm opens to a wider-sense
interpretation where achieving the global objective implies a cost
reduction at every \ug.

\subsection{The cost functions}\label{ssec:cost_fncts}
As mentioned before, the algorithm proposed hereafter works with any set of
generation/transfer cost functions, as long as they all satisfy some mild
convexity constraints. Specifically, it is required that $\{C_i(\cdot)\}$ and
$\gamma(\cdot)$ are positive valued, monotonically increasing, convex and twice
differentiable. Even though these requirements may seem abstract and distant
from real systems, one should take into account that the cost function
of a common electrical generator (as, e.g., oil, coal, nuclear,\dots) is often modeled as
a quadratic polynomial $C(x)=a+bx+cx^2$, where the coefficients $a$, $b$ and $c$
depend on the generator type, see~\cite{Lenoir} and,
especially~\cite{generator_data}. Such a generation model clearly satisfies our
assumptions and, without available counterexamples, we extended it to the
transfer cost function~$\gamma(\cdot)$.

It is worth commenting here that the assumptions on the cost
functions also allow for a simple way to introduce upper bounds on
the energy generated by the \ugs\ or supported by the transfer
connections. Indeed, one can introduce \emph{soft} constraints by
designing the cost function with a steep rise at the nominal maximum
value (see also the first paragraph of
Section~\ref{sec:num_results}, where we comment on the cost
functions used for simulation). By doing so, the maximum
generated/transferred energies are controlled directly by the cost
functions, without the need for hard constraints (i.e., well
specified inequalities such as $E_i^{(g)} \le
E_{i,\mathrm{max}}^{(g)}$) that would increase the complexity of the
minimization problem in (\ref{eq:start}). Note that this expedient
translates, somehow, to a more flexible system: when needed, a \ug\ 
can produce more energy than the nominal maximum if it is willing to
pay an (significant) extra cost. Indeed, this situation arises in
practical systems when backup generators are activated.

\section{Iterative Distributed Minimization}
\label{sec:minimization}
%
\subsection{Decentralizing the problem}
Problem (\ref{eq:start}) is known to have a unique minimum point
since both the objective function and the constraints are strictly
convex. However, dealing with $M(M-1)$ unknowns can be very
involved. Moreover, a centralized solution would require a control
unit that is aware of all the system characteristics. This fact
implies a considerable amount of data traffic to gather all the
information and can miss some privacy requirements, since \ugs\ may
prefer to keep production costs and quantities private. To avoid
these issues, we propose here a distributed iterative approach that
reaches the minimum cost by decomposing the problem into $M$ local,
reduced-complexity subproblems solved by the \ugs\ with little
information about the rest of the system.

\subsubsection{Identifying local subproblems}
In order to decompose~(\ref{eq:start}) into $M$ \ug\ subproblems, let us rewrite
it in the following equivalent form
\begin{equation}\label{eq:step1}
\begin{aligned}
\mathcal{C}^*=\min_{\{\Es_i\},\{E_{i,j}\}} & \sum_{i=1}^M
C_i\bigl(E_i^{(c)}+\Es_i-\ve_i^T\MA^T\vE_i^{(b)}\bigr)
+ \sum_{i=1}\ve_i^T\MA^T\boldsymbol{\gamma}(\vE_i^{(b)})\\
\text{s.\ to}\;\;
          & E_{i,j}\ge 0, \forall i,j,\\
              & E_i^{(c)}+\Es_i-\ve_i^T\MA^T\vE_i^{(b)}\ge
0, \forall i,\\
          &\Es_i = \ve_i^T\MA\vE_i^{(s)}, \forall i.\\
\end{aligned}
\end{equation}
The idea is that, for each \ug, we first use the new variable
$\Es_i$ to represent the energy sold by \ug-$i$ and only later we
force it to be equal to all the energy bought by other \ugs\ from
\ug-$i$, namely $\Es_i = \ve_i^T\MA\vE_i^{(s)}$, the coupling
constraint.

Due to the convexity properties of the primal problem~(\ref{eq:start}) (or,
equivalently,~(\ref{eq:step1})), one can find the minimum cost by
relaxing the $M$ coupling constraints and solving the dual problem
\begin{equation}\label{eq:dual}
\mathcal{C}^* = \max_{\vlambda} \mathcal{C}(\vlambda)
\end{equation}
where $\mathcal{C}(\vlambda) = \sum_{i=1}^M
\mathcal{C}^{(l)}_i(\vlambda)$ with terms
\begin{equation}\label{eq:min_local}
\begin{aligned}
\mathcal{C}^{(l)}_i(\vlambda)=\min_{\Es_i,\vE_i^{(b)}} &
\mathcal{C}_i(\Es_i,\vE_i^{(b)},\vlambda)\\
\text{s.\ to}\;\; & \Es_i\ge 0, E_{j,i}\ge0, \forall j\\
              & E_i^{(c)}+\Es_i-\ve_i^T\MA^T\vE_i^{(b)}\ge0.
\end{aligned}
\end{equation}
In the last definition, which is a local minimization subproblem given the
parameters $\vlambda$, we introduced
\begin{equation}\label{eq:lagrangian}
\mathcal{C}_i(\Es_i,\vE_i^{(b)},\vlambda) = C_i\bigl(E_i^{(c)}+\Es_i-\ve_i^T\MA^T\vE_i^{(b)}\bigr)
{}+ \ve_i^T\MA^T\boldsymbol{\gamma}(\vE_i^{(b)}) +
\ve_i^T\MA^T\diag\{\vlambda\}\vE_i^{(b)}- \lambda_i\Es_i,
\end{equation}
that is the contribution of \ug-$i$ to the Lagrangian function relative
to~(\ref{eq:step1}). The parameter vector
$\vlambda={\begin{bmatrix}\lambda_1&\cdots&\lambda_M\end{bmatrix}}^T$ gathers
all the Lagrange multipliers $\lambda_i$ corresponding to the coupling
constraints $\Es_i=\ve_i^T\MA\vE_i^{(s)}$, respectively and for all
$i=1,\dots,M$.

\subsubsection{Iterative dual problem solution}
To solve the dual problem~(\ref{eq:dual}), we resort to the iterative
subgradient method~\cite[Chapter 8]{Bertsekas_Convex}, which basically finds a
sequence $\{\vlambda[k]\}$ that converges to the optimal point of the dual
problem~(\ref{eq:dual}), namely $\vlambda^*=\arg\max_{\vlambda}\mathcal{C}(\vlambda)$.
More specifically, for each point $\vlambda[k]$, each \ug\ minimizes its
contribution to the Lagrangian function by solving the local
subproblem~(\ref{eq:min_local}) and determining the minimum point
$(\Es_i[k],\vE_i^{(b)}[k])=(\Es_i(\vlambda[k]),\vE_i^{(b)}(\vlambda[k]))$. Then,
the Lagrange multipliers are updated according to
\begin{equation}\label{eq:clearing}
\vlambda[k+1] = \vlambda[k] + \alpha[k]\begin{bmatrix}\ve_1^T\MA\vE_1^{(s)}[k]-
\Es_1[k]\\\vdots\\\ve_M^T\MA\vE_M^{(s)}[k]-\Es_M[k]\end{bmatrix},
\end{equation}
where $\alpha[k]$ is a positive step factor. Also, recall
from~(\ref{eq:en_vectors}) that the set of vectors $\{\vE_i^{(s)}\}$
can be readily derived knowing the set $\{\vE_i^{(b)}\}$. Note that
the vector
$\boldsymbol{\varsigma}={[\ve_i^T\MA\vE_i^{(s)}[k]-[\Es_i[k]]}_{M\times
1}$ is a subgradient of the dual concave function
$\mathcal{C}(\vlambda)$ in $\vlambda=\vlambda[k]$, i.e.\
$\mathcal{C}(\vlambda)\le\mathcal{C}(\vlambda[k]) +
\boldsymbol{\varsigma}^T(\vlambda-\vlambda[k]), \forall \vlambda$.
Finally, (\ref{eq:clearing}) also says that $\lambda_i$ can be
updated at \ug-$i$ once the vector $\vE_i^{(s)}[k]$ has been built
with the inputs $E_{i,j}[k]$ collected from the neighboring \ugs.

\subsubsection{Interpretation---Market clearing}

Algorithm~\ref{alg:distributed} summarizes the distributed
minimization procedure. One can readily notice that all necessary
data is computed at the \ugs, with no need for an external,
centralized control unit. The information exchanged by the \ugs\ is
limited to the Lagrange multipliers $\{\lambda_i\}$ and the demanded
energies $\{E_{j,i}\}$, computed at \ug-$i$ and communicated only to
the corresponding \ug-$j$. Both privacy and traffic limitations are
hence satisfied.

\begin{algorithm}
\caption{Distributed approach}\label{alg:distributed}
\begin{algorithmic}
\State \ug-$i$ initialize $\lambda_i[0]$
\Repeat
    \State \ugs\ exchange $\{\lambda_i[k]\}$
    \State \ug-$i$ computes $\Es_i[k]$ and $\vE_i^{(b)}[k]$ by
      solving~(\ref{eq:min_local}) with fixed $\vlambda[k]$
    \State \ug-$i$ informs \ug-$j$, $j\ne i$,  about the energy it is willing
      to buy, namely $E_{j,i}[k]$, at the given price $\lambda_{j}[k]$
    \State with energy requests $E_{i,j}[k]$ from neighboring \ugs, \ug-$i$
      builds $\vE_i^{(s)}[k] \gets {\begin{bmatrix}E_{i,1}[k] & \cdots & E_{i,M}[k]\end{bmatrix}}^T$
    \State \ug-$i$ computes
      $\lambda_i[k+1] \gets \lambda_i[k] + \alpha[k] (\ve_i^T\MA\vE_i^{(s)}[k]-\Es_i[k])$
    \State $k\gets k+1$
\Until{convergence condition is verified}
\end{algorithmic}
\end{algorithm}

As commented before, this algorithm allows for an interesting
interpretation: each Lagrange multiplier $\lambda_i$ may be
understood as the price per energy unit requested by \ug-$i$ to sell
energy to its neighbors. Then, the Lagrangian
function~(\ref{eq:lagrangian}) can be seen as the ``net
expenditure'' (the opposite of the net income) for \ug-$i$: each
\ug\ pays for producing energy, for buying energy and for
transporting the energy it buys. Conversely, the \ug\ is payed for
the energy it sells. By solving problem~(\ref{eq:min_local}) ,
\ug-$i$ is thus maximizing its benefit for some given selling
($\lambda_i[k]$) and buying ($\lambda_j[k], j\ne i$) prices per
energy unit. According to this view, the updating
step~(\ref{eq:clearing}) is \emph{clearing the market}: prices
should be modified until, globally, energy demand matches energy
offer. Note that (\ref{eq:clearing}) is an example of the law of
demand: if the energy offered by \ug-$i$ $\Es_i[k]$ is less than all
the energy demanded by the neighboring \ugs\ from \ug-$i$, that is
$\ve_i^T\MA\vE_i^{(s)}$, then the selling price must increase and
$\lambda_i[k+1]\ge\lambda_i[k]$.

\subsection{The \protect\ug\ subproblem}
\label{ssec:subproblem}
In the previous section we have shown how the cost minimization
problem~(\ref{eq:start}) can be solved by means of successive
iterations between the solution of local
problems~(\ref{eq:min_local}) and the update of the Lagrange
multipliers according to~(\ref{eq:clearing}). We will give now a
closed-form solution to the local subproblem~(\ref{eq:min_local}) to be solved by the generic \ug-$i$.

In order to keep notation as simple as possible, and without loss of generality,
we assume that the Lagrange multipliers $\{\lambda_j\}, j\ne i$ are ordered in
increasing order, i.e.\ $\lambda_{\min}=\lambda_1\le\lambda_2 \le\cdots
\le\lambda_{i-1} \le\lambda_{i+1} \le\cdots\le\lambda_M$. Also, with some abuse
of notation, we fix $\lambda_j=+\infty$ when $a_{j,i}=0$: as far as \ug-$i$ is
concerned, the fact that there is no connection from \ug-$j$ to \ug-$i$ is
equivalent to assume that the price of the energy sold by \ug-$j$ is too high to
be worth buying. Besides, we will make use of the functions $C_i'(\cdot)$ and
$\gamma'(\cdot)$ (the first derivatives of the cost functions $C_i(\cdot)$ and
$\gamma(\cdot)$) and of their inverse functions, respectively $\chi_i(\cdot)$
and $\Gamma(\cdot)$. It is interesting to mention that, in economics, the
derivative of a cost function is called the \emph{marginal cost} (the cost of
increasing infinitesimally the argument). To see this, consider for instance the
generation cost function and assume that we increase production from
$E_i^{(g)}$ to $E_i^{(g)}+\epsilon$, with $\epsilon$ representing a
small amount of energy. Then, the new generation cost can be
approximated as
\begin{equation}\label{eq:marginal}
C_i(E_i^{(g)}+\epsilon)\approx C_i(E_i^{(g)}) + C_i'(E_i^{(g)})\epsilon,
\end{equation}
showing that the cost varies proportionally with $\epsilon$ and that the
coefficient is $C_i'(E_i^{(g)})$. Note that $\epsilon$ can be either positive
(more energy is produced for, e.g., selling purposes) or negative (because,
e.g., some energy is bought from outside). Analogously, $\gamma'(E_{j,i})$ is
the marginal transportation cost.

The solution to the minimization subproblem~(\ref{eq:min_local}) at
\ug-$i$ behaves according to six different cases. Each case is
characterized by a specific relationship between the
generation/transportation marginal costs $C_i'(\cdot)$ and
$\gamma'(\cdot)$ and the unitary selling/buying prices
$\{\lambda_i\}$. We report next the mathematical definition of the
six cases whereas, in Section~\ref{ssec:economical_interpretation},
some additional comments will help in grasping their
electrical/economical meaning.

\begin{case}[\ug-$i$ neither sells nor buys]\label{case:1}
If $\lambda_i\le C_i'(E_i^{(c)})$ and $\lambda_{\min} \ge C_i'(E_i^{(c)}) -
\gamma'(0)$, then \ug-$i$ will decide to remain in a self-contained state and
generate all and only the energy it consumes. Namely,
\begin{align*}
E_i^{(g)}&=E_i^{(c)}, & \Es_i &= 0, & E_{j,i}&=0\quad\forall j\ne i.
\end{align*}
\end{case}

\begin{case}[\ug-$i$ buys but neither generates nor sells]\label{case:2}
Let us assume that $\lambda_i\le C_i'(0)$ and $\lambda_{\min} \ge
\lambda_i-\gamma'(E_i^{(c)})$. Moreover, we can identify a partition
$\{\mathcal{S}_*,\mathcal{S}_0\}$ of $\{j=1,\dots,M: j\ne i\}$ that satisfies the
following assumptions:\begin{itemize}
\item it exists $\eta>0$ such that\begin{itemize}
\item $\eta > \lambda_j - \lambda_i + \gamma'(0)$ for all $j\in\mathcal{S}_*$;
\item $\eta \le \lambda_j - \lambda_i + \gamma'(0)$ for all $j\in\mathcal{S}_0$;
\item $\eta$ is the unique positive solution to
$$
\sum_{j\in\mathcal{S}_*}\Gamma(\eta-\lambda_j+\lambda_i) = E_i^{(c)};
$$
\end{itemize}
\item either $\exists j\in\mathcal{S}_* : \lambda_j \ge \lambda_i
-\gamma'(0)$ or $\forall j\in\mathcal{S}_*, \lambda_j < \lambda_i
-\gamma'(0)$ and $\sum_{j\in\mathcal{S}_*}\Gamma(\lambda_i-\lambda_j)\le
E_i^{(c)}$;
\item for all $j\in\mathcal{S}_*$, one has $\lambda_{j} \le C_i'(0)
-\gamma'(0)$ and
$$
\sum_{j\in\mathcal{S}_*} \Gamma\Bigl(C_i'(0)-\lambda_j\Bigr) \le E_i^{(c)}.
$$
\end{itemize}
Then, \ug-$i$ buys all and only the energy it consumes, i.e.\ it neither
generates nor sells any energy. More specifically
\begin{align*}
E_{j,i} &= \Gamma(\eta + \lambda_i - \lambda_j)\quad \forall j\in\mathcal{S}_*, &
E_{j,i} &= 0\quad \forall j\in\mathcal{S}_0,\\
\Es_i &= 0 & E_i^{(g)}&=0.
\end{align*}
\end{case}

\begin{case}[\ug-$i$ generates and buys but does not sell]\label{case:3}
Let us assume that $\lambda_i < C_i'(E_i^{(c)})$ while $\lambda_{\min} >
\max\{C_i'(0),\lambda_i\}-\gamma'(E_i^{(c)})$ and $\lambda_{\min} <
C_i'(E_i^{(c)})-\gamma'(0)$. Moreover, we can identify a partition
$\{\mathcal{S}_*,\mathcal{S}_0\}$ of $\{j=1,\dots,M: j\ne i\}$ that satisfies the
following assumptions:\begin{itemize}
\item it exists $\eta>0$ such that\begin{itemize}
\item $\eta > \lambda_j - \lambda_i + \gamma'(0)$ for all $j\in\mathcal{S}_*$;
\item $\eta \le \lambda_j - \lambda_i + \gamma'(0)$ for all $j\in\mathcal{S}_0$;
\item $\eta$ is the unique positive solution to
$$
C_i'\biggl(E_i^{(c)}-\sum_{j\in\mathcal{S}_*}\Gamma(\eta+\lambda_i+\lambda_j)\biggr)
= \eta + \lambda_i
$$
\end{itemize}
\item either $\exists j\in\mathcal{S}_* : \lambda_j \ge \lambda_i
-\gamma'(0)$ or $\forall j\in\mathcal{S}_*, \lambda_j < \lambda_i
-\gamma'(0)$ and $\lambda_i\le C_i'\biggl(E_i^{(c)}-
\sum_{j\in\mathcal{S}_*}\Gamma(\lambda_i-\lambda_j)\biggr)$;
\item either $\exists j\in\mathcal{S}_* : \lambda_{j} > C_i'(0)
-\gamma'(0)$ or $\forall j\in\mathcal{S}_*, \lambda_{j} \le C_i'(0)
-\gamma'(0)$ and
$$
\sum_{j\in\mathcal{S}_*} \Gamma\Bigl(C_i'(0)-\lambda_j\Bigr) > E_i^{(c)}.
$$
\end{itemize}
Then, \ug-$i$ does not sell any energy. Furthermore, it buys some energy to
supplement the local generator and feed all the loads. The exact amounts are as
follows:
\begin{align*}
E_{j,i} &= \Gamma(\eta + \lambda_i - \lambda_j)\quad \forall j\in\mathcal{S}_*, &
E_{j,i} &= 0\quad \forall j\in\mathcal{S}_0,\\
\Es_i &= 0 & E_i^{(g)}&=\chi_i(\eta+\lambda_i).
\end{align*}
\end{case}

\begin{case}[\ug-$i$ generates and sells but does not buy]\label{case:4}
If $\lambda_i>C_i'(E_i^{(c)})$ and $\lambda_{\min}\ge \lambda_i-\gamma'(0)$,
then \ug-$i$ does not buy any energy. Conversely, it generates all the energy it
needs plus some extra energy for the market. More specifically,
\begin{align*}
E_i^{(g)}&=\chi_i(\lambda_i), & \Es_i &= E_i^{(g)}-E_i^{(c)},  & E_{j,i}&=0\quad\forall j\ne i.
\end{align*}
\end{case}

\begin{case}[\ug-$i$ sells and buys but does not generate]\label{case:5}
Assume that $\lambda_i\le C_i'(0)$ and $\lambda_{\min} < \lambda_i -
\gamma'(0)$. Also, let
$\mathcal{S}_*=\{j:\lambda_j<\lambda_i-\gamma'(0)\}$. Note that
$\mathcal{S}_*\ne\emptyset$ since, at least,
$\lambda_{\min}\in\mathcal{S}_*$. Then, \ug-$i$ does not generate
any energy: it buys all the energy it consumes, plus some extra
energy for the market, from all \ugs\ in the set $\mathcal{S}_*$. The
exact amounts are as follows:
\begin{align*}
E_{j,i} &= \Gamma(\lambda_i - \lambda_j)\quad \forall j\in\mathcal{S}_*, &
E_{j,i} &= 0\quad \forall j\notin\mathcal{S}_*,\\
\Es_i &= \sum_{j\in\mathcal{S}_*}E_{j,i}
-E_i^{(c)}, & E_i^{(g)}&=0.
\end{align*}
\end{case}

\begin{case}[\ug-$i$ sells, buys and generates]\label{case:6}
Assume that $\lambda_i>C_i'(0)$, $\lambda_{\min} < \lambda_i-\gamma'(0)$ and
$$
\lambda_i > C_i'\Bigl(E_i^{(c)}-\sum_{j\in\mathcal{S}_*}
\Gamma(\lambda_i-\lambda_j)\Bigr),
$$
where we introduced the set
$\mathcal{S}_*=\{j:\lambda_j<\lambda_i-\gamma'(0)\}$. Then, the
local generator is activated but \ug-$i$ also buys energy from all
\ugs\ in $\mathcal{S}_*$. After feeding all local loads with
$E_i^{(c)}$, some extra energy is left for selling in the market:
\begin{align*}
E_{j,i} &= \Gamma(\lambda_i - \lambda_j)\quad \forall j\in\mathcal{S}_*, &
E_{j,i} &= 0\quad \forall j\notin\mathcal{S}_*,\\
\Es_i &= E_i^{(g)}+\sum_{j\in\mathcal{S}_*}E_{j,i}
-E_i^{(c)}, & E_i^{(g)}&=\chi_i(\lambda_i).
\end{align*}
\end{case}

\begin{IEEEproof}
The proof of these results is a cumbersome convex optimization
exercise. From the Karush-Kuhn-Tucker conditions associated
to~(\ref{eq:min_local}), one must suppose all the different cases
above and realize that the corresponding assumptions are necessary
for each given case. Furthermore, one can also derive the exact
values of all energy flows. Once all cases have been considered, a
careful inspection shows that the derived necessary conditions form
a partition of the hyperplane $\{(\lambda_i, {\{\lambda_j\}}_{j\ne
i})\}$. Hence, the condition are also sufficient, along with
necessary, and the proof is concluded.
All the details are given in the appendix.
\end{IEEEproof}

As a final remark, note that we are assuming a positive load at the
\ugs, i.e.\ $E_i^{(c)}>0, \forall i$. The case where
\mbox{$E_i^{(c)}=0$} may be handled analogously and brings to
similar results. More specifically, Cases~\ref{case:4},~\ref{case:5}
and~\ref{case:6} extend directly because of continuity. Conversely,
Cases~\ref{case:2} and~\ref{case:3} disappear and expand the domain
of Case~\ref{case:1} to $\lambda_i\le C_i'(0)$ and
$\lambda_{\min}\ge\lambda_i-\gamma'(0)$.

\begin{figure}
\centering \begin{subfigure}{.45\columnwidth}\centering
\begin{tikzpicture}
\end{tikzpicture}
\subcaption{Fully connected\label{fig:fullcnctd}}
\end{subfigure}\quad \begin{subfigure}{.45\columnwidth}\centering
\begin{tikzpicture}
\end{tikzpicture}
\subcaption{Ring\label{fig:ring}}
\end{subfigure}\\[2ex]
\begin{subfigure}{.9\columnwidth}\centering
\begin{tikzpicture}
\end{tikzpicture}
\subcaption{Line\label{fig:line}}
\end{subfigure}
\caption{Considered connection topologies.} \label{fig:topologies}
\end{figure}
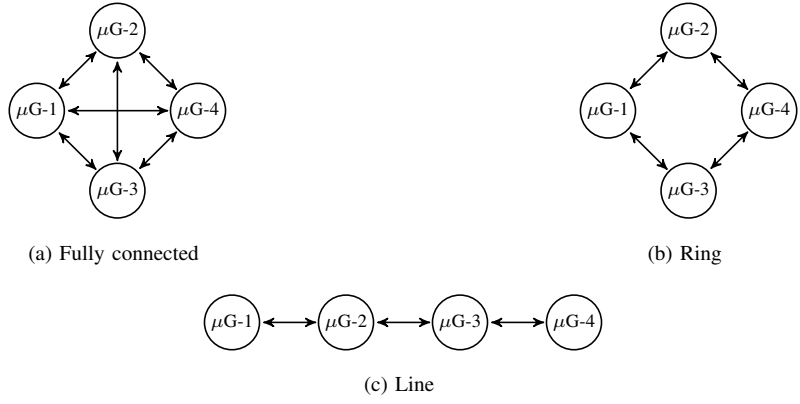

\subsection{Interpretation and summary}
\label{ssec:economical_interpretation}

It is interesting to note that the optimal solution provided in the
previous section has a straightforward interpretation in economical
terms. Recall that, at \ug-$i$, the Lagrange multipliers can be
interpreted as the unitary selling price ($\lambda_i$) and the
unitary buying prices from the other \ugs\ ($\{\lambda_j\}, j\ne i$).
Moreover, the derivatives $C_i'(\cdot)$ and $\gamma'(\cdot)$ are the
marginal generation cost and the marginal transportation cost,
respectively, that is the linear variation on the cost due to an
infinitesimal variation of the generated or transported energy,
respectively.

Bearing this in mind, let us focus on Case~\ref{case:1}. By means
of~(\ref{eq:lagrangian}) and (\ref{eq:marginal}), one readily realizes that \ug-$i$
is not interested in selling energy since the selling price $\lambda_i$ is lower
than the marginal production cost $C_i'(E_i^{(g)})$. Indeed, the income
$\lambda_i\Es_i$ will be lower than the extra production cost, namely
$C_i(E_i^{(c)}+\Es_i)-C_i(E_i^{(c)})>C_i'(E_i^{(c)})\Es_i$ (the last inequality
is due to the convexity of $C_i(\cdot)$). Similarly, buying is not profitable
either since the minimum energy price $\lambda_{\min}$ is larger than the
marginal benefit\footnote{When buying energy, the production cost reduces but
the transportation cost increases. The marginal benefit with respect to
$E_i^{(g)}=E_i^{(c)}$---the \ug\ generates all and only the energy it
consumes---is thus $C_i'(E_i^{(c)})-\gamma'(0)$.} $C_i'(E_i^{(c)})-\gamma'(0)$.
The conditions for Case~\ref{case:1} are hence justified. Analogous
considerations hold for the other cases.

Another interesting point is that microgrids are always willing to
trade since their local cost without trading, i.e. $C_i(E_i^{(c)}) +
\gamma(0)$, will always be higher than their ``net
expenditure''  $\mathcal{C}^{(l)}_i(\vlambda^*)$ in (\ref{eq:lagrangian}), where
$\vlambda^{*}$ stands for the optimal point of
(\ref{eq:dual}). For the sake of brevity, we prove this result for
Case~\ref{case:6} only, although the same reasoning holds true for
the rest of cases. In Case~\ref{case:6}, the ``net expenditure''
reads
$$
  \mathcal{C}^{(l)}_i(\vlambda^*) =
  C_{i}\big(E_i^{(c)}+\Es_i-\sum_{j\in\mathcal{S}_{*}}E_{j,i}\big)
  + \sum_{j\in\mathcal{S}_{*}} \gamma(E_{j,i}) + \sum_{j\in\mathcal{S}_{*}} \lambda_j^*
  E_{j,i} - \lambda_i^* \Es_i
$$
where $(\Es_i,\{E_{j,i}\})=(\Es_i(\vlambda^*),\{E_{j,i}(\vlambda^*)\})$ is now
the minimum point of~(\ref{eq:min_local}) for $\vlambda=\vlambda^*$. Next, by
means of the results of Case~\ref{case:6}, one
has $\lambda_i^* = C_i'(E_i^{(g)})$ and $\lambda_j^* = \lambda_i^* -
\gamma'(E_{j,i})$ for all $j\in \mathcal{S}_{*}$. Then, the
equation above can be rewritten as follows:
$$
  \mathcal{C}^{(l)}_i(\vlambda^*) = C_{i}\big(E_i^{(c)}+E_0\big) - C_{i}'\big(E_i^{(c)}+E_0\big)
  E_0  + \sum_{j\in\mathcal{S}_{*}} \big(\gamma(E_{j,i}) -
  \gamma'(E_{j,i})E_{j,i}\big),
$$
with $E_0 \triangleq \Es_i-\sum_{j\in\mathcal{S}_{*}}E_{j,i}$.
Finally, since the cost functions are monotonically increasing and
convex, it turns out that $C_{i}(E_i^{(c)}+E_0) -
C_{i}'(E_i^{(c)}+E_0)E_0 < C_{i}(E_i^{(c)})$ and $\gamma(E_{j,i}) -
  \gamma'(E_{j,i})E_{j,i} < \gamma(0)=0$, which leads to the desired
  result.

\section{Numerical Results}
\label{sec:num_results}
As for the numerical results, we have considered a system composed
of four microgrids and, for simplicity, we have assumed the same
generation cost function at all microgrids. More specifically, we
have considered the U12 generator in \cite{generator_data}, whose
quadratic cost function $C(x)=a+bx+cx^2$ has coefficients
$a=86.3852$~\$, $b=56.5640$~\$/MWh and $c=0.3284$~\$/(MWh)$^2$.
Besides, as motivated in Section~\ref{ssec:cost_fncts}, the original
cost function has been multiplied by the function $1+{(0.9\cdot
x/E_{\max}^{(g)})}^{30}$ in order to include a soft constraint that
accounts for the maximum energy generation $E_{\max}^{(g)} =
10$~MWh. Regarding the transfer cost function, we have modeled it as
the cubic polynomial\footnote{\label{footnote:transferfunction}
Note that this particular choice is arbitrary and this function may
depend on physical parameters and the business model. However,
similar results are expected with other setups.}
$\gamma(x)=x+x^3$, with $x$ in MWh and $\gamma(x)$ in US dollars.
Numerical results are given for the three topologies shown in
Fig.~\ref{fig:topologies}: fully connected, ring and line.

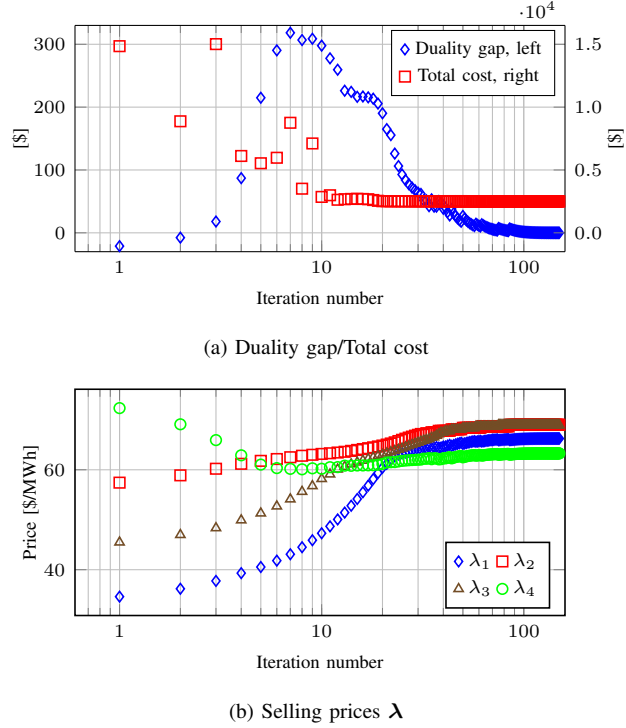
\begin{figure}[!t]
\centering
\begin{subfigure}[b]{.9\columnwidth}
\centering
\begin{tikzpicture}
\end{tikzpicture}
\subcaption{Duality gap/Total cost\label{fig:dualitygap}}
\end{subfigure}\\[2ex]
\begin{subfigure}[b]{.9\columnwidth}
\centering
\begin{tikzpicture}
\end{tikzpicture}
\subcaption{Selling prices $\vlambda$\label{fig:lambdas}}
\end{subfigure}
\caption{Algorithm convergence for $\vE^{(c)} = [8, 11, 11, 6]^T$ MWh:
(\subref{fig:dualitygap}) shows the evolution of the
duality gap (left y-axis) and the evolution of the total cost (right
y-axis), while  (\subref{fig:lambdas}) shows the evolution of the
selling prices $\vlambda$.
\label{fig:Duality_and_lambdas}}
\end{figure}

Fig.~\ref{fig:Duality_and_lambdas} assesses the convergence speed of the
distributed minimization algorithm. The curves refer to a fully connected
system where the microgrid loads are $\vE^{(c)} = [8, 11, 11, 6]^T$ MWh.
First, Fig.~\ref{fig:dualitygap} shows the duality gap and the
total cost of the system as a function of the iteration number. As
we can observe, the algorithm converges (the duality gap is almost null) after a reasonable number of
iterations. More interestingly, Fig.~\ref{fig:lambdas} shows
the evolution of the selling prices. Note that the relationship between
prices after convergence reflects the one between local loads: the more energy the
\ug\ consumes locally, the higher its selling price is. While this is a natural
consequence of the generation cost function being strictly increasing (if the
local load is low, the \ug\ can generate extra energy for selling purposes at a lower cost), we can
also see it as a manifestation of the law of demand. Indeed, the microgrid
energy demand may be seen as composed by two terms, an internal one
corresponding to the local loads and an external one from the other microgrids.
The selling price will hence increase with the resulting total demand. Some more
insights about how and how fast the algorithm converges can be found in \cite{Matamoros2012}.
There, only two microgrids are considered: such a simple
case allows for a centralized closed-form solution and its direct comparison
with the distributed approach.

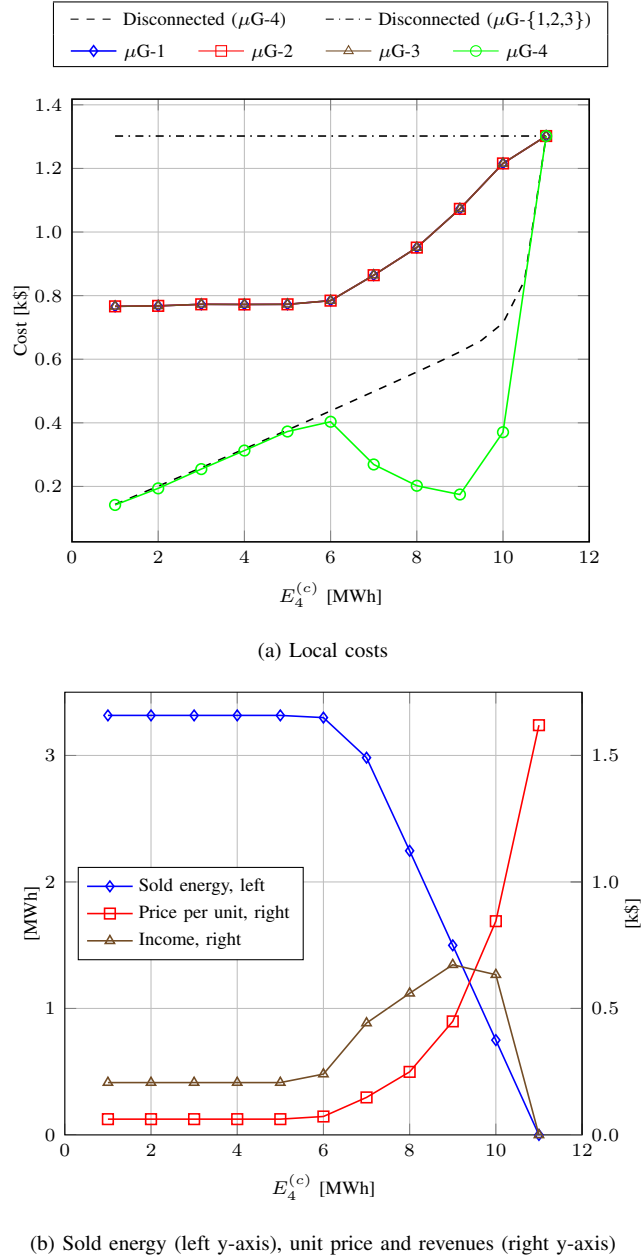
\begin{figure}[!t]
\centering
\begin{subfigure}{0.95\columnwidth}
\centering
\begin{tikzpicture}
\end{tikzpicture}
\caption{Local costs}\label{fig:fully_costs}
\end{subfigure}\\[2ex]
\begin{subfigure}{0.95\columnwidth}
\centering
\begin{tikzpicture}
\end{tikzpicture}
\caption{Sold energy (left y-axis), unit price and revenues (right y-axis)}
\label{fig:impactmu4}
\end{subfigure}
\caption{Fully-connected topology: local costs (\subref{fig:fully_costs})
and \ug-4 metrics (\subref{fig:impactmu4}) for
$\vE^{(c)} = [11, 11, 11, E_4^{(c)}]^T$ MWh.}\label{fig:fully}
\end{figure}

\begin{figure}
\centering
\begin{tikzpicture}
\end{tikzpicture}
\caption{Costs for the ring topology ($\vE^{(c)} = [11, 11, 11, E_4^{(c)}]^T$
MWh).\label{fig:ring_costs}}
\end{figure}
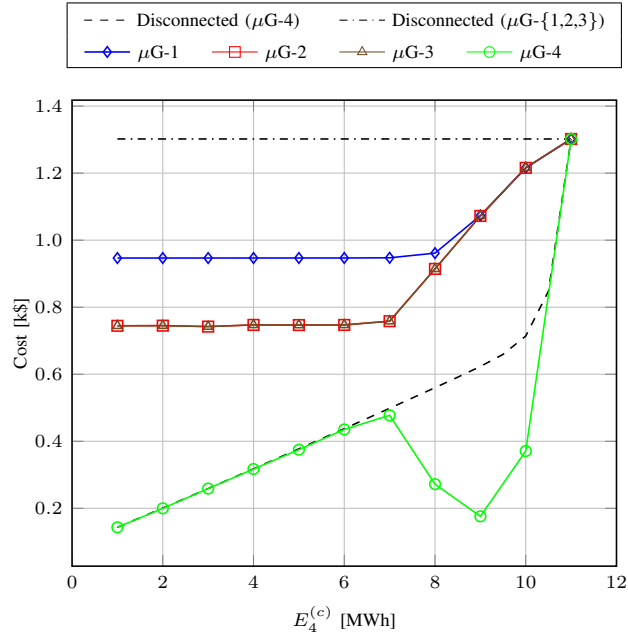

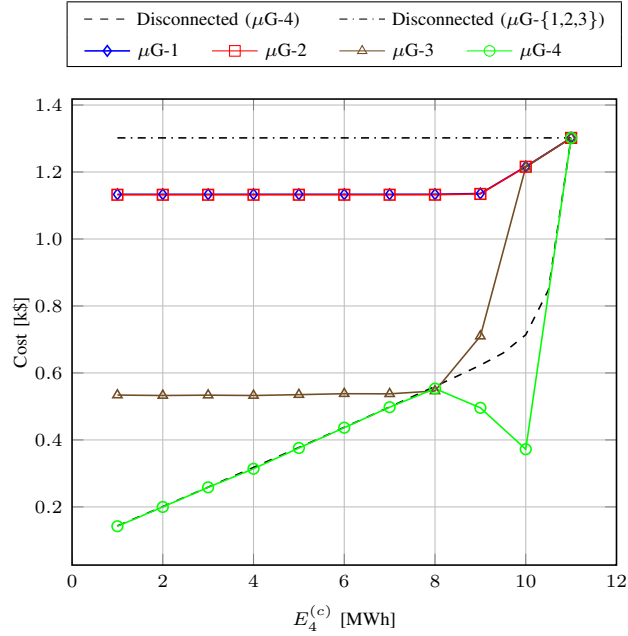
\begin{figure}
\centering
\begin{tikzpicture}
\end{tikzpicture}
\caption{Costs for the line topology ($\vE^{(c)} = [11, 11, 11, E_4^{(c)}]^T$
MWh).\label{fig:line_costs}}
\end{figure}

Next, in order to get some more insight into the evolution of prices
and energy flows, we consider a scenario where all local loads are
held constant at 11~MWh (just above $E_{\max}^{(g)}$), except for
\ug-4, whose load varies from~1 to~11~MWh. In
Figures~\ref{fig:fully_costs},~\ref{fig:ring_costs}
and~\ref{fig:line_costs}, for the four microgrids, we report the
local cost after convergence
$\mathcal{C}_i^{(l)}(\boldsymbol{\lambda}^*)$, that is the minimum
``net expenditure''~(\ref{eq:min_local}), with
$\boldsymbol{\lambda}^*$ the maximum point of (\ref{eq:dual}). For
benchmarking purposes, we have also depicted the costs at each
microgrid in the disconnected case (i.e.\ when no trading is
performed, the dashed lines). As shown in
Section~\ref{ssec:economical_interpretation}, when
using~(\ref{eq:lagrangian}) as the local cost, we can observe that
optimal trading always brings some benefit (cost reduction) to all
microgrids.

Let us now focus on the fully connected topology of
Fig.~\ref{fig:fully}. In Fig.~\ref{fig:fully_costs} we see that the
cost attained by \ug-4 after trading initially follows the cost of
the disconnected microgrid. It is only when the local load grows
above 6~MWh that the gain becomes noticeable, reaches its maximum
for $E_4^{(c)}\approx9$~MWh and then decreases again until it
becomes null at $E_4^{(c)} = 11$~MWh. There, all \ugs\ have the same
internal demand and, for symmetry reasons, there is no energy
exchange. The gain, indeed, is a result of the energy sold by \ug-4
to the other microgrids, whose amount, unit price and corresponding
income is depicted in Fig.~\ref{fig:impactmu4}. At first sight, it
may be disconcerting to see that the negligible gain obtained by
\ug-4 for $E_4^{(c)}<6$~MWh is the result of selling a large amount
of energy at a very low price (both almost constant for
$E_4^{(c)}<6$~MWh). For $E_4^{(c)}>6$~MWh, however, the \ug\ sells
less energy but the unit price increases fast enough to improve the
gain (the ``Income'' curve): the best trade-off between the amount
of energy sold by the \ug\ and its unit price is reached, as said
before, for $E_4^{(c)}\approx9$~MWh. For larger values of
$E_4^{(c)}$, the high selling price cannot compensate for the
decrease of sold energy and the income goes to zero.

To understand why this happens, let us consider \ug-4. After convergence, \ug-4 is
defined by Case~\ref{case:4}---generates and sells. Then,
particularizing~(\ref{eq:lagrangian}), the local costs at \ug-4 is
\begin{equation}\label{eq:cost4}
\mathcal{C}_4=C(E_4^{(c)}+\Es_4)-\lambda_4\Es_4,
\end{equation}
where we used the fact that we are assuming $C_i(\cdot)=C(\cdot)$ for
$i=1,2,3,4$. The optimal $\lambda_4=\lambda_4^*$ is given by the marginal cost, that is
\begin{equation}\label{eq:price4-4}
\lambda_4^* = C'(E_4^{(c)}+\Es_4).
\end{equation}
Now, consider (\ref{eq:cost4}), (\ref{eq:price4-4}) and recall that
the cost function $C(E_4^{(g)})$ is the one with label
``Disconnected (\ug-4)'' in Fig.~\ref{fig:fully_costs}: for
$E_4^{(g)}<9$~MWh approximately, the cost function is almost linear.
Thus, the unit price (\ref{eq:price4-4}) is nearly constant and so
is the cost (\ref{eq:cost4}) as a function of $\Es_4$ for all
$E_4^{(c)}$ and $\Es_4$ such that $E_4^{(g)}=E_4^{(c)}+\Es_4<9$~MWh.
In other words, the \ug\ cost $\mathcal{C}_4$ only depends on the
consumed energy and not on the sold one since the income from
selling some extra energy is canceled out by the extra cost needed
to generate it. All these considerations are reflected by the curves
in Figures ~\ref{fig:fully_costs} and~\ref{fig:impactmu4} for
$E_4^{(c)}<6$~MWh or, equivalently, $E_4^{(g)}<9$~MWh, since the
total energy sold by \ug-4 in this regime is approximately 3~MWh.

Now, without loss of generality, let us focus on \ug-1. Note that, for symmetry
reasons, \ug-1, \ug-2 and \ug-3 are all buying the same amount of energy from
\ug-4 (namely, $E_{4,i}=\Es_4/3$) and they are not exchanging energy to one another.
Intuitively, \ug-1 falls within Case~\ref{case:3}---generates and buys---and its
local cost is
$$
\mathcal{C}_1=C(E_1^{(c)}-\Es_4/3)+\gamma(\Es_4/3)+\lambda_4\Es_4/3.
$$
The optimal price $\lambda_4=\lambda_4^*$ is given by (\ref{eq:price4-4}) but,
from \ug-1 perspective, can also be rewritten as
\begin{equation}\label{eq:price4-1}
\lambda_4^* = C'(E_1^{(c)}-\Es_4/3)-\gamma'(\Es_4/3).
\end{equation}
Let us neglect, for the moment, the transfer cost function $\gamma(\cdot)$.
Also, recall that $C(\cdot)$ is again the one depicted in
Fig.~\ref{fig:fully_costs} with label ``Disconnected (\ug-4)''. Since
$E_1^{(c)}=11$~MWh, the production cost $C(E_1^{(c)}-\Es_4/3)$ takes values
above the curve elbow for all $E_{4,1}=\Es_4/3<1$~MWh, approximately. In this
regime, the generation cost is very high and \ug-1 is certainly trying to buy
energy to reduce its expenditure. For $E_{4,1}=\Es_4/3>1$~MWh, however, the
generation cost takes values below the curve elbow and has an almost linear
behavior (see also comments above). Then, it is not worth to buy more energy,
since its price will cancel out the generation savings. The convex nature of
the transfer cost function $\gamma(\cdot)$ accentuates this trend. These
considerations explain why $\Es_4=3\cdot E_{4,1}\approx 3$~MWh for all
$E_4^{(c)}<6$~MWh.

For $E_4^{(c)}>6$~MWh, the four \ugs\ work in the nonlinear part of
the generation cost function. Matching the energy price at both the
selling side (\ref{eq:price4-4}) and the buying side
(\ref{eq:price4-1}) results in the trade-off of
Fig.~\ref{fig:impactmu4}, which turns out to be quite fruitful for
\ug-4. This fact compensates somehow for the little benefit (with
respect to the benefits of the other \ugs) experienced by \ug-4 at
low values of $E_4^{(c)}$.

Finally, let us recall that the purpose of the original
problem~(\ref{eq:start}) is to minimize the total cost of the system
and not to maximize local benefits. This, together with the fact
that we do not allow \ugs\ to cheat, explains why \ug-4 always sells
at a unit price given by the marginal cost and does not look for
extra gains.

Even though the general ideas discussed above about the cost behaviors apply to
all systems, connection topologies (see Fig.~\ref{fig:topologies}) other than the full connected one present
some specific characteristics. For instance, Fig.~\ref{fig:ring_costs} report the
local cost for the four \ugs\ in ring topology. We may see that \ug-2 and \ug-3
get some extra benefit from acting as intermediaries between \ug-4 and \ug-1.
To be more precise, after the trading
process, the local solutions at \ug-$\{2,3\}$ will fall within Case
\ref{case:6}: energy is bought not only to satisfy
internal needs but also to be resold to \ug-1.  By doing so,
\ug-$\{2,3\}$ can reduce their local cost substantially.

In Fig.~\ref{fig:line} we depict another situation where topology
heavily impacts on the costs. In this case, microgrids are connected
by means of a line topology with \ug-4 at the end of the
line. Therefore, \ug-3 can be
regarded as the bottleneck of the system, since all the energy that goes
from \ug-4 to \ug-\{1,2\} has to pass through it inevitably. As it
can be observed in the figure, this situation benefits \ug-3.

\section{Conclusions}
\label{sec:conclusions}

In this paper, we have addressed a problem in which several
microgrids interact by exchanging energy in order to minimize the
global operation cost, while still satisfying their local demands.
In this context, we have proposed an iterative distributed algorithm
that is scalable in the number of microgrids and keeps local cost
functions and local consumption private. More specifically, each
algorithm iteration consists of a local minimization step followed
by a market clearing process. During the first step, each microgrid
computes its local energy bid and reveals it to its potential
sellers. Next, during the market clearing process, energy prices are
adjusted according to the law of demand. As for the local
optimization problem, it has been shown to have a closed form
expression which lends itself to an economical interpretation. In
particular, we have shown that no matter the local demand that a
microgrid will be always willing to start the trading process since,
eventually, its ``net expenditure'' will be lower than its local
cost when operating on its own. Finally, numerical results have
confirmed that the algorithm converges after a reasonable number of
iterations and there certainly is a gain over nonconnected \ugs\
which strongly depends on the energy demands and network topology.

\appendices
\section{Solution to the Microgrid Problem}
This appendix proves the solution to the \ug\ problem given in
Section~\ref{ssec:subproblem}. As explained before, we will suppose that the
\ug\ best option is to operate in one of the six different states defined
according to whether the \ug\ is (or is not) selling, buying and generating any
energy, as described by the six cases of Section~\ref{ssec:subproblem}. By doing
so, one can compute all the energy values of interest and identify what
constraints the prices $\{\lambda_i,{\{\lambda_j\}}_{j\ne i}\}$ must satisfy for
the considered case to be feasible. After all cases have been considered, the
six sets of necessary conditions should form a partition of the
$(\lambda_i,{\{\lambda_j\}}_{j\ne i})$ hyperplane. Indeed, this fact implies
that each set of conditions is sufficient, along with necessary, for the
corresponding \ug\ state and that the computed energy values are those
minimizing the local cost function (\ref{eq:lagrangian}) for given prices
$\{\lambda_i,{\{\lambda_j\}}_{j\ne i}\}$.

Before delving into the different cases, some common preliminaries are needed.
For the local problem (\ref{eq:min_local}) the Lagrangian function is as follows:
\begin{multline*}
\mathcal{L}=C_i\bigl(E_i^{(c)}+\Es_i
- \ve_i^T\MA^T\vE_i^{(b)}\bigr) + \ve_i^T\MA^T\boldsymbol{\gamma}(\vE_i^{(b)}) +
  \ve_i^T\MA^T\boldsymbol{\Lambda}\vE_i^{(b)}\\{} - \lambda_i\Es_i
-\eta\Es_i-\boldsymbol{\mu}^T\vE_i^{(b)} - \omega\bigl(E_i^{(c)}+\Es_i
- \ve_i^T\MA^T\vE_i^{(b)}\bigr),
\end{multline*}
where we have introduced the Lagrange multipliers $\eta, \omega$ and
$\boldsymbol{\mu}={\begin{bmatrix}\mu_1 &\mu_2&\cdots&\mu_M\end{bmatrix}}^T$.
The KKT conditions hence write
\begin{subequations}\label{eq:KKT}
\begin{align}
&\frac{\partial\mathcal{L}}{\partial \Es_i} =
C_i'\bigl(E_i^{(c)}+\Es_i-\ve_i^T\MA^T\vE_i^{(b)} \bigr) - \lambda_i - \eta - \omega=0\label{eq:KKT1}\\
& \frac{\partial\mathcal{L}}{\partial\vE_i^{(b)}} = -C_i'\bigl(E_i^{(c)}+\Es_i -
\ve_i^T\MA^T\vE_i^{(b)}\bigr)\MA\ve_i \nonumber\\&\hspace{2cm} {}+
\diag\{\MA\ve_i\}\boldsymbol{\gamma}'(\vE_i^{(b)})
+ \mathbf{\Lambda}\MA\ve_i - \boldsymbol{\mu}
+ \omega\MA^T\ve_i=\mathbf{0}\label{eq:KKT2}\\
&\Es_i\ge 0,\quad \eta\ge 0,\quad \eta\Es_i=0,\label{eq:KKT3}\\
&E_{j,i}\ge 0,\quad \mu_j\ge 0,\quad \mu_jE_{j,i}=0,\quad\forall j=1,\dots,M,\label{eq:KKT4}\\
&E_i^{(c)}+\Es_i - \ve_i^T\MA^T\vE_i^{(b)}\ge0,\quad
\omega\ge 0,\quad \omega\bigl(E_i^{(c)}+\Es_i-\ve_i^T\MA^T\vE_i^{(b)})=0.\label{eq:KKT5}
\end{align}
\end{subequations}
By recalling the definition of $\MA$ in Section~\ref{sec:system_model}, the elements of the gradient in
(\ref{eq:KKT2}) can be written in a much simpler form, namely
\begin{align*}
& C_i'\bigl(E_i^{(c)}+\Es_i - \ve_i^T\MA^T\vE_i^{(b)}\bigr)
- \gamma'(E_{j,i}) - \lambda_j + \mu_j - \omega= 0, &&\forall j:a_{j,i}=1,\\
&\qquad \mu_j=0,\quad E_{j,i}=0,	&&\forall j:a_{j,i}=0.
\end{align*}
The derivation of the six possible solutions given in
Section~\ref{ssec:subproblem} is based on the analysis of the KKT conditions
above, as explained hereafter.

\subsection{Proof of Case~\protect\ref{case:1}}
Let us suppose that the solution of the minimization problem tells us that the
\ug\ neither sells nor buys any energy, that is $\Es_i=0$,
$\vE_i^{(b)}=\mathbf{0}$ and $E_i^{(g)}=E_i^{(c)}$. If this was the case, then
the KKT conditions~(\ref{eq:KKT}) would write
\begin{align*}
&C_i'\bigl(E_i^{(c)}\bigr)-\lambda_i-\eta = 0,\\
&C_i'\bigl(E_i^{(c)}\bigr)-\gamma'(0)-\lambda_j + \mu_j = 0, &&\forall j\ne i,\\
&\eta\ge 0, \omega=0 \text{ and } \mu_j\ge 0, &&\forall j\ne i.
\end{align*}
The first condition implies
$$
\eta=C_i'\bigl(E_i^{(c)}\bigr)-\lambda_i\ge0 \Rightarrow \lambda_i\le
C_i'\bigl(E_i^{(c)}\bigr),
$$
while the second condition yields
$$
\mu_j=\gamma'(0)+\lambda_j-C_i'\bigl(E_i^{(c)}\bigr)\ge 0 \Rightarrow
\lambda_j\ge C_i'\bigl(E_i^{(c)}\bigr)-\gamma'(0),
$$
which are the two necessary conditions corresponding to Case~\ref{case:1}.

\subsection{Proof of Case~\protect\ref{case:2}}
We now look for the necessary conditions for Case~\ref{case:2}, that is
\ug-$i$ sells no energy (i.e.\ $\Es_i=0$) and buys energy from at least another
\ug\ (i.e.\ $\exists j=1,\dots,M, j\ne i:E_{j,i}>0$). Besides, \ug-$i$
generates no energy and $E_i^{(g)}=E_i^{(c)}-\ve_i^T\MA^T\vE_i^{(b)}=0$. The KKT
conditions~(\ref{eq:KKT}) simplify to
\begin{subequations}
\begin{align}
&C_i'(0)-\lambda_i-\eta-\omega = 0,\label{eq:KKT1_2}\\
&C_i'(0)-\gamma'(E_{j,i})-\lambda_j-\omega = 0 &&\text{and } \mu_j=0,
	&&\forall j\in\mathcal{S}_*,\label{eq:KKT2_2}\\
&C_i'(0)-\gamma'(0)-\lambda_j+\mu_j-\omega = 0 &&\text{and } \mu_j\ge 0,
	&&\forall j\in\mathcal{S}_0,\label{eq:KKT3_2}\\
&\eta\ge 0\text{ and }\omega\ge 0,\label{eq:KKT4_2}
\end{align}
where we have introduced the sets
\begin{align*}
\mathcal{S}_* &= \{j=1,\dots,M:j\ne i \text{ and } E_{j,i}>0\},\\
\mathcal{S}_0 &= \{j=1,\dots,M:j\ne i \text{ and } E_{j,i}=0\}.
\end{align*}
\end{subequations}

The first necessary condition $\lambda_i\le C_i'(0)$ is a
straightforward consequence of (\ref{eq:KKT1_2}) and (\ref{eq:KKT4_2}). Next,
for all $j\in\mathcal{S}_*$, (\ref{eq:KKT1_2}) and (\ref{eq:KKT2_2}) imply
\begin{equation}\label{eq:5.2a}
\lambda_i+\eta-\lambda_j-\gamma'(E_{j,i})=0,
\end{equation}
which means that
\begin{equation}\label{eq:5.2b}
\lambda_j<\lambda_i+\eta-\gamma'(0),\quad\forall j\in\mathcal{S}_*.
\end{equation}
Similarly, because of (\ref{eq:KKT1_2}) and (\ref{eq:KKT3_2}), one has
$$
\lambda_j\ge\lambda_i+\eta-\gamma'(0),\quad\forall j\in\mathcal{S}_0.
$$
By comparing the last two inequalities, one sees that $\lambda_j<\lambda_k$ for
all $j\in\mathcal{S}_*,k\in\mathcal{S}_0$, meaning that
$\lambda_{\min}=\min_{j\ne i} \lambda_j$ is certainly part of the set
$\mathcal{S}_*$ when this solution is correct.

From (\ref{eq:5.2a}), and taking (\ref{eq:5.2b}) into account, we can infer
that
$$
E_{j,i}=\Gamma(\eta-\lambda_j+\lambda_i),
$$
where $\Gamma(\cdot)$ is the inverse of $\gamma'(\cdot)$, which exists because
of the continuity and convexity assumptions on the cost function $\gamma(\cdot)$.
Since we are supposing that the optimal working point satisfies
$E_i^{(c)}-\ve_i^T\MA^T\vE_i^{(b)}=0$, we see that $\eta$ shall satisfy the
equality
\begin{equation}\label{eq:5.2c}
\sum_{j\in\mathcal{S}_*}\Gamma(\eta-\lambda_j+\lambda_i) = E_i^{(c)}.
\end{equation}
Given that $\Gamma(\cdot)$ is an increasing function, and recalling
(\ref{eq:5.2b}), it can be easily shown that equation (\ref{eq:5.2c}) in
the variable $\eta$ has a unique solution, whose value allows us to compute the
energies $E_{j,i}$ bought from neighbor \ug-$j$, $j\in\mathcal{S}_*$ (compare
with the statement of Case~\ref{case:2}).

To derive the other necessary conditions for this case, let us focus on
(\ref{eq:5.2c}). Since $\Gamma(\cdot)$ is a non-negative function
and $\lambda_j=\lambda_{\min}\Rightarrow j\in\mathcal{S}_*$, it follows:
$$
\Gamma(\eta-\lambda_{\min}+\lambda_i) \le E_i^{(c)}.
$$
Recalling that $\Gamma(\cdot)$ is the inverse of $\gamma'(\cdot)$, this implies
$0\le\eta\le\gamma'\bigl(E_i^{(c)}\bigr)+\lambda_{\min}-\lambda_i$
and, hence, $\lambda_{\min}\ge\lambda_i-
\gamma'\bigl(E_i^{(c)}\bigr)$ is a necessary condition for Case~\ref{case:2}.

Now, let $\mathcal{S}_{\dagger}=\{j=1,\dots,M:j\ne i\text{ and } \lambda_j <
\lambda_i - \gamma'(0)\}$, which is a subset of $\mathcal{S}_*$ because of
(\ref{eq:5.2b}). Then, by comparison with (\ref{eq:5.2c}),
$$
E_i^{(c)} \ge \sum_{j\in\mathcal{S}_{\dagger}}\Gamma(\eta-\lambda_j+\lambda_i)
  \ge \sum_{j\in\mathcal{S}_{\dagger}}\Gamma(\lambda_i-\lambda_j),
$$
where the second inequality is a consequence of $\Gamma(\cdot)$ being an
increasing function. Also, we are letting $\eta\to 0$, which is possible since
$\lambda_i-\lambda_j>\gamma'(0)$ for the considered $j$. In particular,
if $\mathcal{S}_{\dagger}=\mathcal{S}_*$ then
$\sum_{j\in\mathcal{S}_*}\Gamma(\lambda_i-\lambda_j) \le E_i^{(c)}$, as reported
in Case~\ref{case:2}. We will see later that this condition is
important to identify the boundary between the solution regions of
Case~\ref{case:2} and Case~\ref{case:5}.

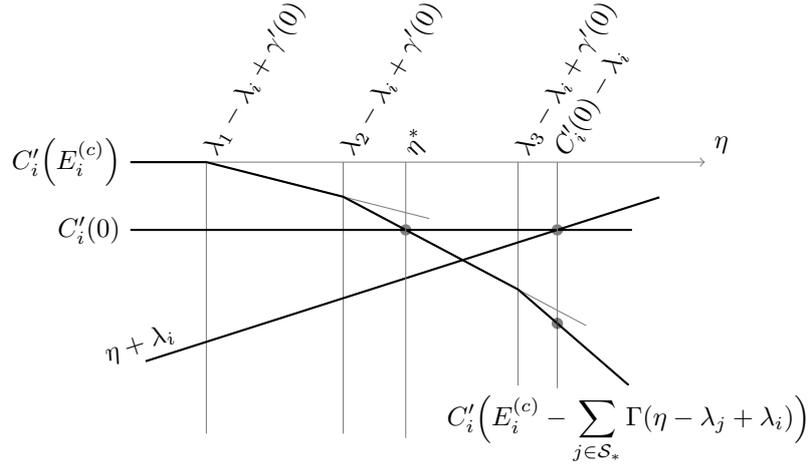
\begin{figure}
\centering
\begin{tikzpicture}
\end{tikzpicture}
\caption{Graphical representation of inequality (\ref{eq:5.2e}) and of the solution $\eta=\eta^*$ to
$C_i'(0)=C_i'\bigl(E_i^{(c)}-\sum_{j\in\mathcal{S}_*} \Gamma(\eta-\lambda_j
+\lambda_i)\bigr)$. Without loss of
generality, we assume here that $\lambda_1=\lambda_{\min}\le\lambda_2\le
\dots\le\lambda_{i-1}\le\lambda_{i+1}\le\dots\le\lambda_M$. Note that, as $\eta$
increases, a new mode is activated each time a point
$\lambda_j-\lambda_i+\gamma'(0)$ is crossed.\label{fig:help2}}
\end{figure}

Consider now (\ref{eq:KKT1_2}) and (\ref{eq:5.2b}), which give
$$
0\le\omega=C_i'(0)-\lambda_i-\eta<C_i'(0)-\lambda_j-\gamma'(0)
$$
and, in turn, a new necessary condition: 
\begin{equation}\label{eq:5.2d}
\lambda_j<C_i'(0)-\gamma'(0),\forall j\in\mathcal{S}_*.
\end{equation}
Moreover, (\ref{eq:KKT1_2}) and (\ref{eq:5.2c}) yield
\begin{equation}\label{eq:5.2e}
C_i'(0)=C_i'\biggl(E_i^{(c)}-\sum_{j\in\mathcal{S}_*} \Gamma(\eta-\lambda_j
+\lambda_i)\biggr)\ge\eta+\lambda_i.
\end{equation}
Substituting into (\ref{eq:5.2c}), we can write the last necessary
condition of Case~\ref{case:2}, namely
\begin{equation}\label{eq:5.2f}
\sum_{j\in\mathcal{S}_*}\Gamma\bigl(C_i'(0)-\lambda_j\bigr) \ge E_i^{(c)}.
\end{equation}
Note that the left-hand side has a meaning according to (\ref{eq:5.2d}). The last
inequality may also be deduced from the graphical representation in
Fig.~\ref{fig:help2}.

\subsection{Proof of Case~\protect\ref{case:3}}
We suppose again a solution where $\Es_i=0$ and $E_{j,i}>0$ for at least one
$j=1,\dots,M,j\ne i$. As opposed to the previous case, however, we also suppose
that $E_i^{(c)}-\ve_i^T\MA^T\vE_i^{(b)}>0$, i.e.\ $\mu$G-$i$ produces some
energy. With this solution, the Lagrange multiplier $\omega$ is zero and the KKT
conditions become
\begin{subequations}
\begin{align}
&C_i'(E_i^{(c)}-\ve_i^T\MA^T\vE_i^{(b)})-\lambda_i-\eta = 0,\label{eq:KKT1_3}\\
&C_i'(E_i^{(c)}-\ve_i^T\MA^T\vE_i^{(b)})-\gamma'(E_{j,i})-\lambda_j = 0 &&\text{and } \mu_j=0,
	&&\forall j\in\mathcal{S}_*,\label{eq:KKT2_3}\\
&C_i'(E_i^{(c)}-\ve_i^T\MA^T\vE_i^{(b)})-\gamma'(0)-\lambda_j+\mu_j = 0 &&\text{and } \mu_j\ge 0,
	&&\forall j\in\mathcal{S}_0,\label{eq:KKT3_3}\\
&\eta\ge 0,\nonumber
\end{align}
where, again,
\begin{align*}
\mathcal{S}_* &= \{j=1,\dots,M:j\ne i \text{ and } E_{j,i}>0\},\\
\mathcal{S}_0 &= \{j=1,\dots,M:j\ne i \text{ and } E_{j,i}=0\}.
\end{align*}
\end{subequations}

\begin{figure}
\centering
\begin{tikzpicture}
\end{tikzpicture}
\caption{Graphical representation of inequality (\ref{eq:5.3d}) and of the solution $\eta=\eta^*$ to
$\eta+\lambda_i=C_i'\bigl(E_i^{(c)}-\sum_{j\in\mathcal{S}_*} \Gamma(\eta-\lambda_j
+\lambda_i)\bigr)$. Without loss of
generality, we assume here that $\lambda_1=\lambda_{\min}\le\lambda_2\le
\dots\le\lambda_{i-1}\le\lambda_{i+1}\le\dots\le\lambda_M$. Note that, as $\eta$
increases, a new mode is activated each time a point
$\lambda_j-\lambda_i+\gamma'(0)$ is crossed.\label{fig:help3}}
\end{figure}
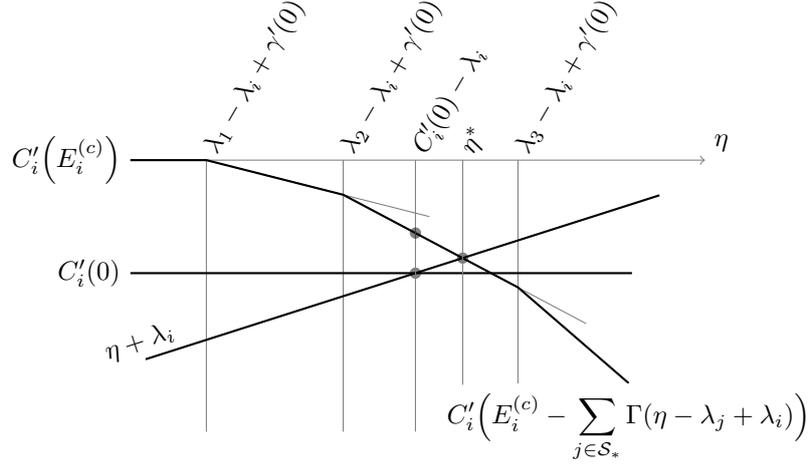

Condition (\ref{eq:KKT1_3}) directly gives the first requirement for the
case in hand, namely $\lambda_i < C_i'(E_i^{(c)})$. Next, by
combining (\ref{eq:KKT1_3}) and (\ref{eq:KKT3_3}), one has
$\lambda_i-\lambda_j+\eta-\gamma'(0)+\mu_j =0$ and, hence, 
\begin{equation}\label{eq:5.3z}
\lambda_j\ge\lambda_i+\eta-\gamma'(0),\quad\forall j\in\mathcal{S}_0,
\end{equation}
since $\mu_j\ge0$. Similarly, (\ref{eq:KKT1_3}) and (\ref{eq:KKT2_3}) yield
$\lambda_i-\lambda_j+\eta-\gamma'(E_{j,i}) =0$, which implies
\begin{equation}\label{eq:5.3a}
\lambda_j < \lambda_i + \eta - \gamma'(0),\quad\forall j\in\mathcal{S}_*,
\end{equation}
since $\gamma'(\cdot)$ is an increasing function, and
\begin{equation}\label{eq:5.3b}
E_{j,i} = \Gamma(\eta+\lambda_i-\lambda_j),
\end{equation}
where $\Gamma(\cdot)$ is the inverse of $\gamma'(\cdot)$. Inequality
(\ref{eq:5.3a}) guarantees that $E_{j,i}$ is positive. Following, by injecting
(\ref{eq:5.3b}) into (\ref{eq:KKT1_3}), it turns out that $\eta$ must satisfy
\begin{equation}\label{eq:5.3c}
C_i'\biggl(E_i^{(c)}-\sum_{j\in\mathcal{S}_*}\Gamma(\eta+\lambda_i-\lambda_j)\biggr) =
\eta+\lambda_i.
\end{equation}
It is straightforward to show that the last equation in $\eta$ admits a unique
solution (see also the graphical representation in Fig.~\ref{fig:help3}) that
allows us to compute the energies bought from neighbor microgrids according to
(\ref{eq:5.3b}), as stated by Case~\ref{case:3}. Also, knowing $\eta$ and
recalling that $\chi_i(\cdot)$ is the inverse function of $C_i'(\cdot)$,
(\ref{eq:KKT1_3}) allows us to compute the generated energy
$$
E_i^{(g)}=E_i^{(c)}-\sum_{j\in\mathcal{S}_*}E_{j,i}=\chi_i(\eta+\lambda_i).
$$

Combining (\ref{eq:5.3c}) with (\ref{eq:5.3a}), we get $\lambda_j+\gamma'(0) <
C_i'(E_i^{(c)})$ for all $\lambda_j, j\in\mathcal{S}_*$ and in particular for
$\lambda_{\min}=\min_{j}\lambda_j$ (see the statement of Case~\ref{case:3}).
We are particularly interested in $\lambda_{\min}$ since the corresponding
microgrid will certainly belong to the set $\mathcal{S}_*$ (and possibly be its
only element) if this case is the solution to the minimization problem. This can
be deduced by comparing (\ref{eq:5.3z}) and (\ref{eq:5.3a}).

Another necessary condition for this solution can be derived by noting that
$E_i^{(c)}-\sum_{j\in\mathcal{S}_*}\Gamma(\eta+\lambda_i-\lambda_j)>0$ implies
$E_i^{(c)}-\Gamma(\eta+\lambda_i-\lambda_{\min})>0$ and, hence,
$\eta<\gamma'(E_i^{(c)})-\lambda_i+\lambda_{\min}$. This bound requires
$\lambda_{\min}>\lambda_i-\gamma'(E_i^{(c)})$ and, together with
(\ref{eq:5.3c}), $\lambda_{\min}>C_i'(0)-\gamma'(E_i^{(c)})$. 

From (\ref{eq:5.3c}) one further has $\eta>C_i'(0)-\lambda_i$. Then, again
because of $\Gamma(\cdot)$ being increasing in $\eta$,
\begin{align}
C_i'(0) &< C_i'\biggl(E_i^{(c)}-\sum_{j\in\mathcal{S}_*}
\Gamma(\eta+\lambda_i-\lambda_j)\biggr)\nonumber\\
&< C_i'\biggl(E_i^{(c)}-\sum_{j\in\mathcal{S}_{\ddagger}}
\Gamma(\eta+\lambda_i-\lambda_j)\biggr)\nonumber\\
&< C_i'\biggl(E_i^{(c)}-\sum_{j\in\mathcal{S}_{\ddagger}}
\Gamma\bigl(C_i'(0)-\lambda_j\bigr)\biggr),\label{eq:5.3d}
\end{align}
and, equivalently, $\sum_{j\in\mathcal{S}_{\ddagger}}
\Gamma\bigl(C_i'(0)-\lambda_j\bigr)<E_i^{(c)}$, where we have introduced
$\mathcal{S}_{\ddagger}=\{j=1,\dots,M:j\ne i\text{ and }\lambda_j<C_i'(0)-
\gamma'(0)\}\subseteq\mathcal{S}_*$ (see also Fig.~\ref{fig:help3}). In
particular, if $\mathcal{S}_{\ddagger}=\mathcal{S}_*$ then
$\sum_{j\in\mathcal{S}_*} \Gamma\bigl(C_i'(0)-\lambda_j\bigr)<E_i^{(c)}$, which
represents the complementary set of (\ref{eq:5.2f}) in
Case~\ref{case:2}.

Finally, let $\mathcal{S}_{\dagger}=\{j=1,\dots,M:j\ne i\text{ and }\lambda_j <
\lambda_i-\gamma'(0)\}$. Then,
\begin{align*}
\lambda_i&\le\eta+\lambda_i=C_i'\biggl(E_i^{(c)}-\sum_{j\in\mathcal{S}_*}
\Gamma(\eta+\lambda_i-\lambda_j)\biggr)\\
&\le C_i'\biggl(E_i^{(c)}-\sum_{j\in\mathcal{S}_{\dagger}}
\Gamma(\eta+\lambda_i-\lambda_j)\biggr)\\
&\le C_i'\biggl(E_i^{(c)}-\sum_{j\in\mathcal{S}_*}
\Gamma(\lambda_i-\lambda_j)\biggr),
\end{align*}
which is equivalent to $\chi_i(\lambda_i) + \sum_{j\in\mathcal{S}_\dagger}
\Gamma(\lambda_i-\lambda_j) \le E_i^{(c)}$. In particular, note that $\chi_i(\lambda_i) + \sum_{j\in\mathcal{S}_*}
\Gamma(\lambda_i-\lambda_j) \le E_i^{(c)}$ when $\mathcal{S}_{\dagger} =
\mathcal{S}_*$ and compare it with Case~\ref{case:6} below.

\subsection{Proof of Case~\protect\ref{case:4}}
The next potential solution provides $E_{j,i}=0$ for all $j\ne i$ (\ug-$i$ does
not buy any energy), $\Es_i>0$ and,
obviously, $E_i^{(g)}=E_i^{(c)}+\Es_i>0$ (the \ug\ sells and generates some
energy). According to the KKT conditions, we also have
$\eta=\omega=0$ and
\begin{align*}
&C_i'\bigl(E_i^{(c)}+\Es_i\bigr)-\lambda_i=0,\\
&C_i'\bigl(E_i^{(c)}+\Es_i\bigr)-\gamma'(0)-\lambda_j+\mu_j=0,&&\forall j\ne i.
\end{align*}

The first KKT condition can be satisfied only if
$\lambda_i>C_i'(E_i^{(c)})$ (the first necessary condition for the current case). Moreover,
combining both conditions, we get $\mu_j=\lambda_j+\gamma'(0)-\lambda_j$, which
requires $\lambda_j\ge\lambda_i-\gamma'(0)$ (the second necessary
condition for the current case) in order to have $\mu_j\ge0$.

The value taken by $\Es_i$ is obtained by inverting $C_i'(\cdot)$ in the first
condition above, namely $\Es_i=\chi_i(\lambda_i)-E_i^{(c)}$. This also means
$E_i^{(g)}=\chi_i(\lambda_i)$.

\subsection{Proof of Case~\protect\ref{case:5}}
We now analyze potential solutions of the type $\Es_i>0$, $E_{j,i}>0$ for at least
one index $j\ne i$ and $E_i^{(c)}+\Es_i-\ve_i^T\MA^T\vE_i^{(b)}=0$. In other
terms, the microgrid does not generate any energy but buys more than consumes
and sells the surplus. In this case, the KKT conditions write
\begin{subequations}
\begin{align}
&C_i'(0)-\lambda_i-\omega = 0,\label{eq:KKT1_5}\\
&C_i'(0)-\gamma'(E_{j,i})-\lambda_j-\omega = 0 &&\text{and } \mu_j=0,
	&&\forall j\in\mathcal{S}_*,\label{eq:KKT2_5}\\
&C_i'(0)-\gamma'(0)-\lambda_j+\mu_j-\omega = 0 &&\text{and } \mu_j\ge 0,
	&&\forall j\in\mathcal{S}_0,\label{eq:KKT3_5}\\
&\eta= 0\text{ and }\omega\ge 0,\nonumber
\end{align}
where we have introduced the sets
\begin{align*}
\mathcal{S}_* &= \{j=1,\dots,M:j\ne i \text{ and } E_{j,i}>0\},\\
\mathcal{S}_0 &= \{j=1,\dots,M:j\ne i \text{ and } E_{j,i}=0\}.
\end{align*}
\end{subequations}

Since $\omega\ge0$, condition (\ref{eq:KKT1_5}) requires
$\lambda_i\le C_i'(0)$, the first necessary condition for the current case. Next, joining
(\ref{eq:KKT1_5}) and (\ref{eq:KKT3_5}), we readily see that
$\mu_j=\lambda_j-\lambda_i+\gamma'(0)$ and
\begin{align}\label{eq:5.5a}
\lambda_j&\ge\lambda_i-\gamma'(0) &&\forall j\in\mathcal{S}_0.
\end{align}

Similarly, (\ref{eq:KKT1_5}) and (\ref{eq:KKT2_5}) imply
$\lambda_i-\lambda_j-\gamma'(E_{j,i})=0$ for all $j\in\mathcal{S}_*$. By using
the function $\Gamma(\cdot)$, inverse of $\gamma'(\cdot)$, the last equation
allows computing the value of the energy bought from $\mu$G-$j$, namely
$E_{j,i}=\Gamma(\lambda_i-\lambda_j)$, as we wanted to show. Note, however, that
this is a meaningful expression only if
$\lambda_j<\lambda_i-\gamma'(0)$ for all $j\in\mathcal{S}_*$. By comparing this
last requirement with (\ref{eq:5.5a}), we see that
$j\in\mathcal{S}_* \Leftrightarrow \lambda_j<\lambda_i - \gamma'(0)$, which is
the definition of the set $\mathcal{S}_*$ as stated by Case~\ref{case:5}.

Finally, the total sold energy is
$$
\Es_i = \ve_i^T\MA^T\vE_i^{(b)} - E_i^{(c)} =
\sum_{j\in\mathcal{S}_*}\Gamma(\lambda_i-\lambda_j)-E_i^{(c)}.
$$
Since we require $\Es_i>0$, it must be
$\sum_{j\in\mathcal{S}_*}\Gamma(\lambda_i-\lambda_j) > E_i^{(c)}$,
our last necessary condition (compare it with Case~\ref{case:2}).

\subsection{Proof of Case~\protect\ref{case:6}}
The last type of solution we have to deal with is characterized by $\Es_i>0$,
$E_{j,i}>0$ for at least one index $j\ne i$ and
$E_i^{(g)}=E_i^{(c)}+\Es_i-\ve_i^T\MA^T\vE_i^{(b)}>0$, i.e.\ \ug-$i$ generates,
buys and sells certain positive amounts of energy. The KKT conditions become
\begin{subequations}
\begin{align}
&C_i'\Bigl(E_i^{(c)}+\Es_i-\ve_i^T\MA^T\vE_i^{(b)}\Bigr)-\lambda_i =
0,\label{eq:KKT1_6}\\
&C_i'\Bigl(E_i^{(c)}+\Es_i-\ve_i^T\MA^T\vE_i^{(b)}\Bigr)-\gamma'(E_{j,i})-\lambda_j = 0 &&\text{and } \mu_j=0,
        &&\forall j\in\mathcal{S}_*,\label{eq:KKT2_6}\\
&C_i'\Bigl(E_i^{(c)}+\Es_i-\ve_i^T\MA^T\vE_i^{(b)}\Bigr)-\gamma'(0)-\lambda_j+\mu_j = 0 &&\text{and } \mu_j\ge 0,
	&&\forall j\in\mathcal{S}_0,\label{eq:KKT3_6}\\
&\eta= 0\text{ and }\omega= 0,\nonumber
\end{align}
where we have introduced the sets
\begin{align*}
\mathcal{S}_* &= \{j=1,\dots,M:j\ne i \text{ and } E_{j,i}>0\},\\
\mathcal{S}_0 &= \{j=1,\dots,M:j\ne i \text{ and } E_{j,i}=0\}.
\end{align*}
\end{subequations}

Condition (\ref{eq:KKT1_6}) directly implies $\lambda_i > C_i'(0)$ (the first
necessary condition for Case~\ref{case:6})
because of the convexity assumptions on $C(\cdot)$. Also, for all
$j\in\mathcal{S}_0$, the Lagrangian multipliers $\mu_j$ are given by
(\ref{eq:KKT1_6}) and (\ref{eq:KKT3_6}), namely
$\mu_j=\lambda_j-\lambda_i+\gamma'(0)$. This value is non-negative only if
\begin{align}\label{eq:5.6a}
\lambda_j&\ge\lambda_i-\gamma'(0) && \forall j\in\mathcal{S}_0.
\end{align}

On the other hand, (\ref{eq:KKT1_6}) and (\ref{eq:KKT2_6}) lead to
$\lambda_i-\lambda_j-\gamma'(E_{j,i})=0$ for all $j\in\mathcal{S}_*$. This
identity allows us to express the quantity
$E_{j,i}=\Gamma(\lambda_i-\lambda_j)$ (the desired result), with $\Gamma(\cdot)$ the inverse of
$\gamma(\cdot)$, provided that $\lambda_j<\lambda_i-\gamma'(0)$ for all
$j\in\mathcal{S}_*$. This requirement, together with (\ref{eq:5.6a}), results
in the definition of $\mathcal{S}_*$ for Case~\ref{case:6}, namely
$j\in\mathcal{S}_* \Leftrightarrow \lambda_j<\lambda_i - \gamma'(0)$.

To conclude, by means of the function $\chi_i(\cdot)$, inverse of $C_i'(\cdot)$,
(\ref{eq:KKT1_6}) is equivalent to $E_i^{(g)}=\chi_i(\lambda_i)$ or, also,
$\Es_i=\chi_i(\lambda_i) -E_i^{(c)} +
\sum_{j\in\mathcal{S}_*}\Gamma(\lambda_i-\lambda_j)$, which are the required
expressions of the total generated and sold energy, respectively. Note that this
is a positive (and meaningful) quantity only if
$\sum_{j\in\mathcal{S}_*}\Gamma(\lambda_i-\lambda_j)+
\chi_i(\lambda_i)>E_i^{(c)}$, the last necessary condition for
Case~\ref{case:6} (compare also with Case~\ref{case:3}).

\subsection{Summary and remarks}
\begin{figure}
\centering
\begin{tikzpicture}
\end{tikzpicture}
\caption{\label{fig:areas}Example of partitioning of the
$(\lambda_i,\{\lambda_j\})$-space according to the solution regions of
Cases~\ref{case:1} to~\ref{case:6} of Section~\ref{ssec:subproblem}. Note that the figure only represents a slice at a
given $(\lambda_i,\lambda_{\min})$-plane: this is enough to identify most of the
boundaries, which only depend on these two parameters. Some boundaries (the
dashed ones), however, depend on all $\{\lambda_j\}$ and cannot be 
represented properly.}
\end{figure}
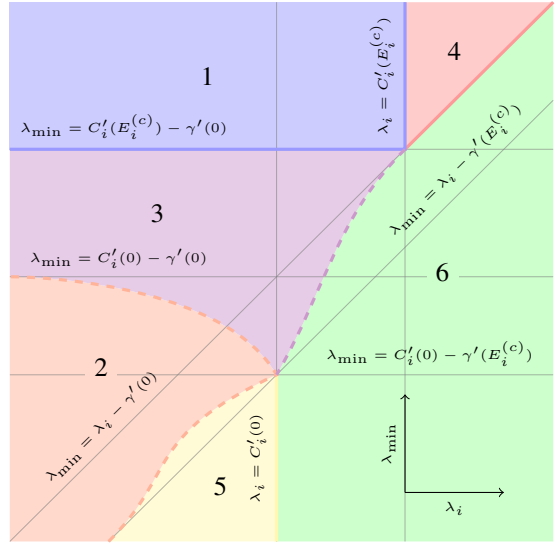

So far, we have analyzed all the cases representing possible operating states
for the microgrid. For each case we have found the corresponding values taken by the
energy flows $E_i^{(g)}$, $\Es_i$ and $E_{j,i}, j\ne i$. Also, each case is
characterized by a set of conditions that are necessary for the feasibility of
the case itself. A close inspection of these
conditions shows that they are mutually exclusive (see also the representation
of Fig.~\ref{fig:areas}). Thus, each set of conditions is also sufficient, together
with necessary, for the respective solution case, meaning that the local minimization
problem is univocally solved.

\end{document}